\documentclass[graybox]{svmult}

\usepackage{mathptmx}       
\usepackage{helvet}         
\usepackage{courier}        
\usepackage{type1cm}        
\usepackage{makeidx}         
\usepackage{graphicx}        
\usepackage{multicol}        
\usepackage[bottom]{footmisc}

\usepackage{cite}
\usepackage{fancyvrb}        

\begin{document}

\title*{Adaptive wavelet-based method for simulation of electronic circuits}
\author{Kai Bittner
\and
Emira Dautbegovic
}
\institute{Kai Bittner\at University of Applied Science of Upper Austria, 4232 Hagenberg,
            Austria,
  \email{kai.bittner@fh-hagenberg.at}
  \and Emira Dautbegovic \at Infineon Technologies AG, 81726 Munich, Germany,
  \email{Emira.Dautbegovic@infineon.com}
}

%
%
\maketitle

\abstract*{In this paper we present an algorithm for analog simulation of
 electronic circuits involving a spline Galerkin method with wavelet-based
 adaptive refinement. Numerical tests show that a first algorithm prototype,
 build within a productively used in-house circuit simulator,
 is completely able to meet and even surpass the accuracy requirements
 and has a performance close to
 classical time-domain simulation methods, with high potential for further
 improvement.
}

\abstract{In this paper we present an algorithm for analog simulation of
 electronic circuits involving a spline Galerkin method with wavelet-based
 adaptive refinement. Numerical tests show that a first algorithm prototype,
 build within a productively used in-house circuit simulator,
 is completely able to meet and even surpass the accuracy requirements
 and has a performance close to
 classical time-domain simulation methods, with high potential for further
 improvement.
}

\section{Introduction}

Wavelet theory emerged during the $20^\mathrm{th}$ century from the
study of Calderon-Zygmund operators in mathematics, the study of the
theory of subband coding in engineering and the study of
renormalization group theory in physics. The common foundation for
the wavelet theory was laid down at the end of the 80's and
beginning of the 90's by work of Daubechies \cite{bittner:D88,bittner:D92}, 
Morlet
and Grossman \cite{bittner:GM84}, Donoho \cite{bittner:D93}, 
Coifman \cite{bittner:CW92},
Meyer \cite{bittner:M93}, Mallat \cite{bittner:M98} and others. Today wavelet-based
algorithms are already in productive use in a broad range of
applications 
\cite{bittner:K94,bittner:MNG+04,bittner:M98,bittner:M93,bittner:MMO+08,
bittner:P03,bittner:WK95,bittner:Y93}, such as
image and signal compression (JPEG2000 standard, FBI fingerprints
database), speech recognition, numerical analysis (solving operator
equations, boundary value problems), stochastics,
smoothing/denoising data, physics (molecular dynamics, geophysics,
turbulence), medicine (heart-rate and ECG analysis, DNA analysis) to
name just a few. Recent approaches
\cite{bittner:BaKnPu08,bittner:CS01,bittner:DCB05,bittner:SGN07,bittner:ZC99} 
to the problem of multirate
envelope simulation indicate that wavelets could also be used to
address the qualitative challenge by a development of novel
wavelet-based circuit simulation techniques capable of an efficient
simulation of mixed analog-digital circuits \cite{bittner:ED_SCEE08}.

The wavelet expansion of a function $f$ is given as
\begin{equation}\label{bittner:wavelet_expansion}
f= \sum_{k\in\mathcal{I}} c_k\, \phi_k +
\sum_{j=0}^\infty \sum_{k\in\Lambda_j}d_{jk}\,\psi_{jk}.
\end{equation}
Here, $j$ refers to a level of resolution, while $k$ describes the
localization in time or space, i.e., $\psi_{jk}$ is essentially
supported in the neighborhood of a point $x_{jk}$. The wavelet
expansion can be seen as coarse scale approximation
$\sum_{k\in\mathcal{I}} c_k\, \phi_k$ by the scaling functions $\phi_k$
complemented by detail information of increasing resolution $j$
in terms of the wavelets $\psi_{jk}$.

In the classical theory wavelets are generated as translation and
dilations of a mother wavelet $\psi$, i.e., $\psi_{jk}(x) =
\psi(2^{-j}x-k)$. However, more general approaches are often used,
e.g., for the construction of wavelets on the interval \cite{bittner:DKU2}
or wavelets for finite element spaces \cite{bittner:NgSt03}. In particular,
non-uniform spline wavelets \cite{bittner:Bit05b} will be used in our
wavelet-based circuit simulation.

Since a wavelet basis consist of an infinite number of wavelets, in
practical computations one has to consider approximations of $f$ by
partial sums of the wavelet expansion
(\ref{bittner:wavelet_expansion}). A simple approach is to fix a
maximal wavelet level $J$ and approximate $f$ by
\begin{equation}\label{bittner:wavelet_linear}
f_J= \sum_{k\in\mathcal{I}} c_k\, \phi_k + \sum_{j=0}^J\sum_{k\in\Lambda_j}
d_{jk}\,\psi_{jk}.
\end{equation}
This approach is called linear approximation, since the
approximation is determined in the linear space of wavelets with
level less or equal $J$. For wavelets of sufficient regularity, one
obtains error estimates of the form
\begin{equation}\label{bittner:approx_linear}
\|f - f_J\|_{L_2} \le C\, 2^{-J s}\,\|f\|_{W_2^s},
\end{equation}
with the Sobolev space $W_2^s$. However, approximation results as
(\ref{bittner:approx_linear}) hold also for other approximation methods,
e.g., for Fourier sums (see \cite{bittner:DeV}).

The real approximation power of wavelets is due to their locality,
which implies that (\ref{bittner:approx_linear}) holds also for small
subintervals. Thus, a piecewise smooth function can be
essentially approximated by some coarse scale approximation with
wavelets added only at non-smooth parts to achieve a required accuracy.
Doing this adaptively for any given signal leads to the notion of
best $n$-term approximation, where the approximation is determined
as linear combination of $n$ arbitrarily chosen wavelets. This results in
an essentially improved approximation for a wide class of
functions, e.g., piecewise smooth function with isolated
singularities. For details about this adaptive, nonlinear
approximation methods we refer to \cite{bittner:DeV,bittner:DVLo}. Usually it is
not obvious which wavelets have to be chosen for optimal
approximation results. In practice optimal wavelet representations can be 
determined by one of the two complementary strategies: coarsening
or refinement.

Coarsening is used if one has already a fine, highly accurate but
expensive approximation, e.g.,  from measurements. The goal is to
throw away as much information as possible, while introducing only a
small error. For a wavelet representation this can be achieved quite
easily by thresholding, which means that wavelets with small
expansion coefficients are removed from the representations.
Inherent stability properties of wavelets ensure that this elimination of 
terms with small
coefficients do not add up to a significant error of the wavelet expansion. 
For wavelets with good localization and
approximation properties, one has many small wavelet coefficients
for piecewise smooth signals with few local singularities (e.g.,
sharp transients), which will result in an essential reduction of
data for the coarsened signal. A disadvantage of coarsening is that
it might be too costly to acquire the fine representation. In
particular, for solving operator equations as in circuit simulation
the reason for using adaptive wavelet techniques is the reduction of
computational cost, which is thwarted by computing a non-adaptive
solution in advance.

In contrast, the strategy of refinement is to start with coarse
approximation and introduce successively more and more degrees of
freedom (e.g., wavelets) in order to improve the approximation.
However, since it is not known in advance, where refinements are
necessary, one has to rely on rough estimates. Therefore it is
reasonable to do the refinement in several steps. This allows to
check the previous steps, while acquiring more information for later
steps. This approach is in particular interesting for iterative
methods, where the approximation is improved in each iteration step
and the number of degrees of freedom can be increased accordingly.

\section{An Adaptive Wavelet Galerkin Method}

We consider circuit equations in the charge/flux oriented
modified nodal analysis (MNA) formulation, which yields a mathematical model
in the form of an initial-value problem of differential-algebraic equations
(DAEs):
\begin{equation}
  \label{bittner:eq_MNA_charge}
  \frac{d}{dt}\vec{q}\big(\vec{x}(t)\big)
        + \vec{f}\big(\vec{x}(t)\big) - \vec{s}(t) = 0.
\end{equation}
Here
\vec{x} is the vector of node potentials and specific branch currents and
\vec{q} is the vector of charges and fluxes. Vector \vec{f} comprises static
contributions, while \vec{s} contains the contributions of independent
sources.

In our adaptive wavelet approach we first discretize the MNA equation
(\ref{bittner:eq_MNA_charge}) in terms of the wavelet basis functions, by
expanding $\vec{x}$ as a linear combination of wavelets
or related functions, i.e., $\vec{x}=\sum_{k=0}^n \vec{c}_k\,\varphi_k$.
For such $\vec{x}$  we integrate  the circuit equations against test functions
 $\theta_\ell$ and obtain the equations
\begin{equation}
\label{bittner:nonlinear}
\int_0^T \Big(\frac{d}{dt}\vec{q}\big(\vec{x}(t)\big)
        + \vec{f}\big(\vec{x}(t)\big) - \vec{s}(t)\Big)\,\theta_\ell\,dt = 0,
\end{equation}
for $\ell=1,\ldots,n$.
Together with the initial conditions $\vec{x}(0)= \vec{x}_0$, we have now
$n+1$ vector valued equations, which determine the coefficients $\vec{c}_i$
provided that the test functions
$\theta_\ell$ are chosen suitably to the basis functions
$\varphi_i$.

Due to the intrinsic properties of wavelets \cite{bittner:ED_SCEE08}
nonlinear wavelet approximation can provide an efficient
representation of functions with steep transients, which often
appear in a mixed analog/digital electronic circuit. However, for an
efficient circuit simulation we have to take into account further
properties of a wavelet system. We consider spline wavelets to be
the optimal choice since spline wavelets are the only wavelets with
an explicit formulation. This permits the fast computation of
function values, derivatives and integrals, which is essential for
the efficient numerical solution of a nonlinear problem as given in
(\ref{bittner:eq_MNA_charge}) (see also \cite{bittner:BiUr04, bittner:DSX00}).
Spline wavelets have been already used for circuit simulation
\cite{bittner:ZhCa99}. However, here we use a completely new approach based
on spline wavelets from \cite{bittner:Bit05b}.

With a good initial guess, Newton's method is known to converge quadratically.
However, a good initial guess is usually not available. In practice we can often
obtain convergence only with a slow converging initial phase of damped Newton steps,
which will mainly contribute to the computational cost of the problem.
On the other
hand, to get a good approximation of the solution of
(\ref{bittner:eq_MNA_charge}), the space $X =\mathrm{span}\{\varphi_k:~k=0,\ldots,n\}$
has to be sufficiently large and the computational cost of each step depends
on $n=\dim X$.
Our approach is to use adaptive wavelet refinement during the Newton iteration,
which leads to an efficient adaptive representation and essentially reduced
computation time.

\section{Interval Splitting Method \label{bittner:split}}

A prototype of the proposed adaptive wavelet algorithm is
implemented within the framework of a productively used circuit
simulator and tested on a variety of circuits. For tests on some
typical RF circuits (amplifier, mixer, oscillator), we were able to
reproduce the results from the transient analysis of the same
circuit simulator up to high accuracy (see \cite{bittner:BiDau10a}). For all
these examples, the wavelet method used a considerable smaller grid
(i.e.\  larger stepsize)
than the transient analysis, while the computation time was higher
but still close to the standard method. This shows that there is a
potential for wavelet methods in circuit simulation, if further
optimization can be achieved.

\begin{figure}
\centering
\includegraphics[width=0.5\textwidth]{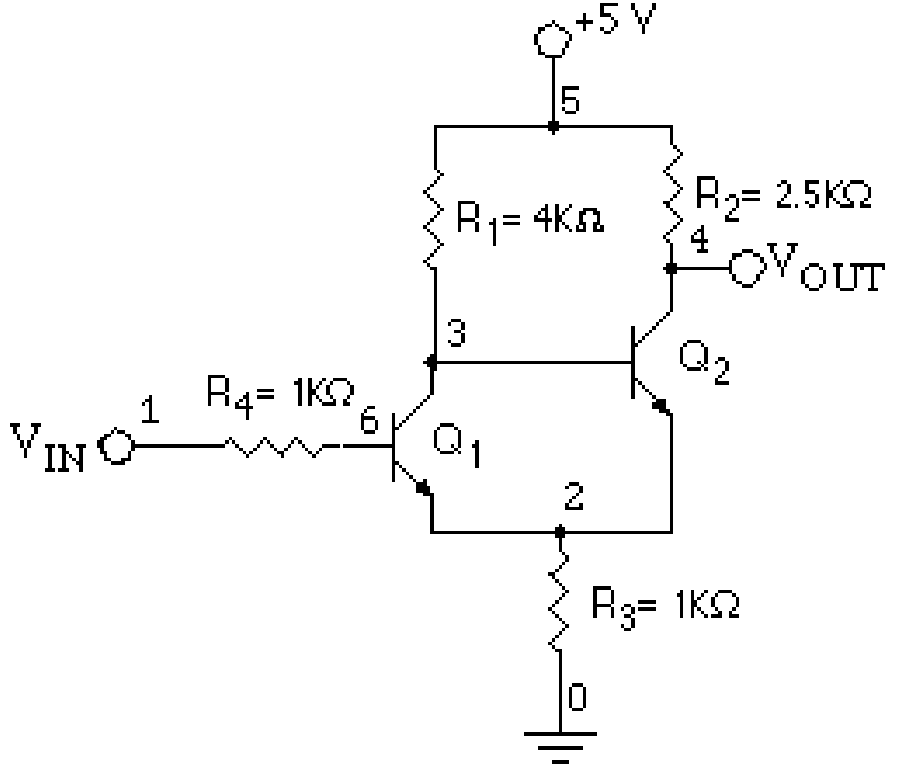}
\caption{Schematics of Schmitt trigger from \cite{bittner:schmitt}}
\label{bittner:schmitt}       
\end{figure}

However, in further tests with a Schmitt trigger circuit
(Fig.~\ref{bittner:schmitt}, \cite{bittner:schmitt}) 
convergence could only be achieved with
a highly accurate initial guess. This is of limited practical value,
since we can usually not provide an initial guess of such quality.
We identified inherit hysteresis of the Schmitt trigger as the main
cause for this problem. In circuits exhibiting hysteresis, certain
input voltages can result in different output, depending on the
previous behaviour of the input signal. With an insufficient initial
guess Newton's method may approach locally the wrong result. This
effect was observed in a Harmonic Balance simulation too, where the
solution is also represented by a basis expansion over an entire
period.

This convergence problem was successfully addressed by a further
improvement of the basic wavelet method based on an interval
splitting mechanism. Basically the wavelet method is applied to a
series of smaller intervals when no convergence is detected. This is
an analogous approach to the reduction of the step size in transient
analysis if no convergence is encountered in the current time step.
In order to preserve continuity, the initial value for each interval
is obtained from the wavelet expansion of the solution on the
previous interval. Furthermore, the interval size is adapted after
each successful step, aiming to keep the problem size for the
wavelet method in a nearly optimal range.

\section{Numerical Tests \label{bittner:numtest}}

The interval splitting method was implemented as an enhancement to
the basic wavelet algorithm and tested on a variety of circuits. For
all examples we have compared the CPU time and the grid size (i.e.,
the number of spline knots or time steps) with the corresponding
results from transient analysis of the underlying circuit simulator.

The error is estimated by comparison with well established
high accuracy transient analysis. The estimate shown in the signal is the
maximal absolute difference over all transient grid points, which gives a good
approximation of the maximal error. That is, if we can obtain a small
error for the wavelet analysis, this proves good agreement with the standard
method. In particular, since we compare the solutions of two independent methods
we have very good evidence that we approximate the solution of underlying
DAE's with the estimated error.

\begin{figure}[htbp]
\centering
\begin{tabular}{c}
\includegraphics[width=0.9\textwidth]{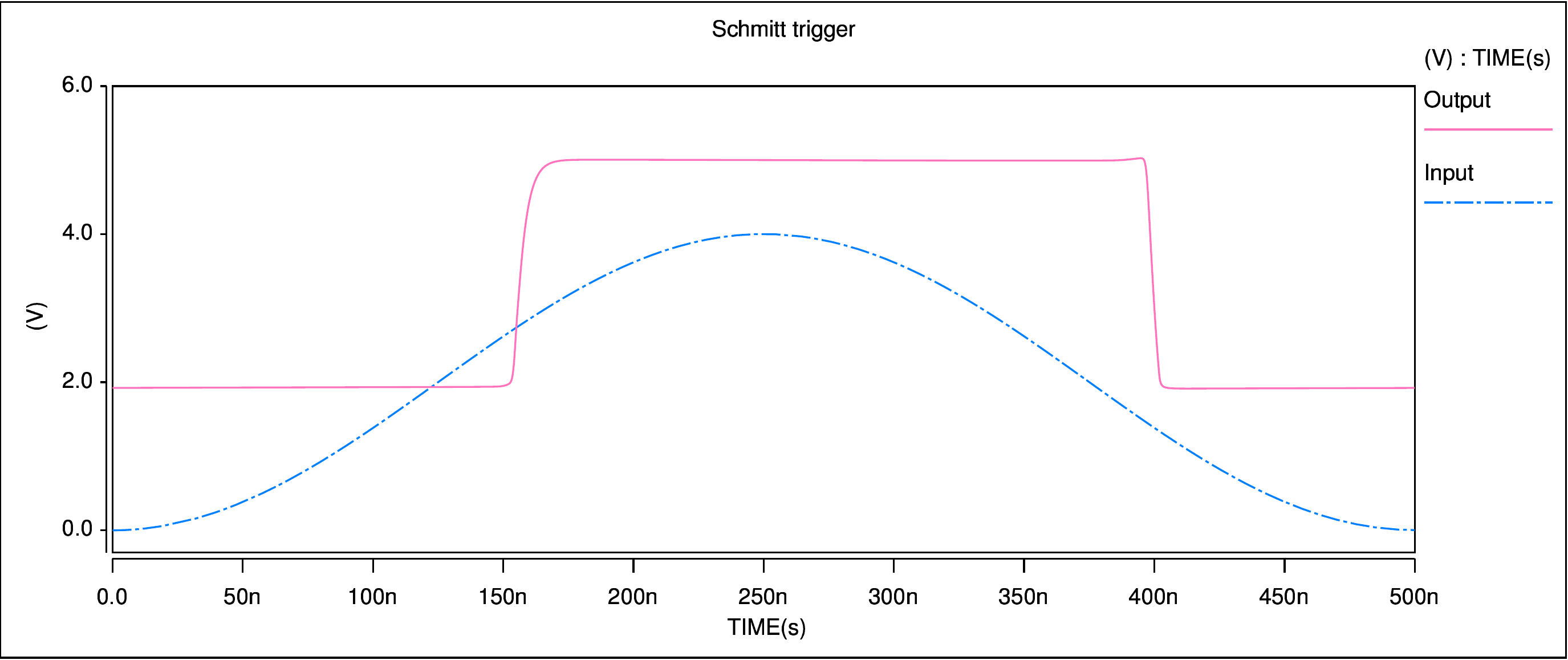}
\end{tabular}
\caption{Input and output signal for the Schmitt trigger.
\label{bittner:Schmitt}}
\end{figure}

{\bf Schmitt trigger.}
The first test circuit is the Schmitt trigger \cite{bittner:schmitt}.
As can be seen in Fig.~\ref{bittner:Schmitt} the output of the Schmitt trigger
circuit signal jumps to a higher level
if the input exceeds an upper threshold and jumps back to low if the input
falls below a lower threshold. However, due to capacitances present in the
transistor model the jumps are slightly smoothed and delayed.

\begin{figure}[htbp]
\centering
\begin{tabular}{cc}
\includegraphics[width=0.48\textwidth]{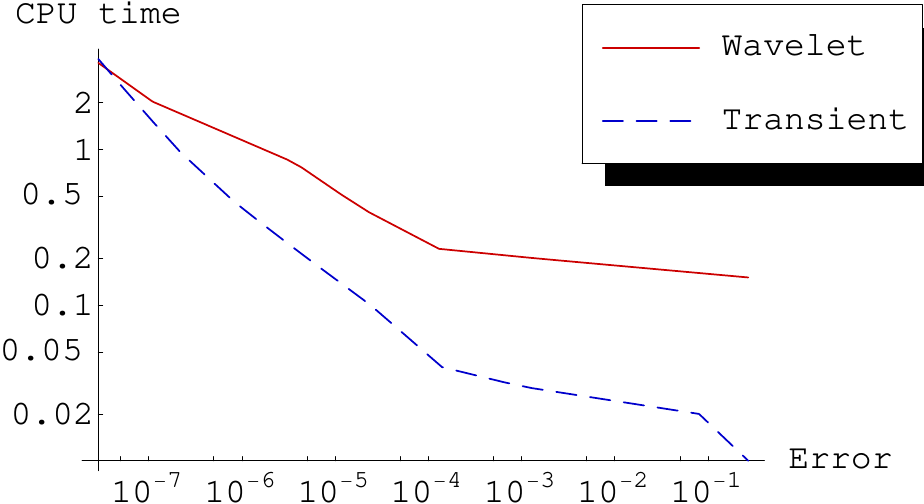}
 &\includegraphics[width=0.48\textwidth]{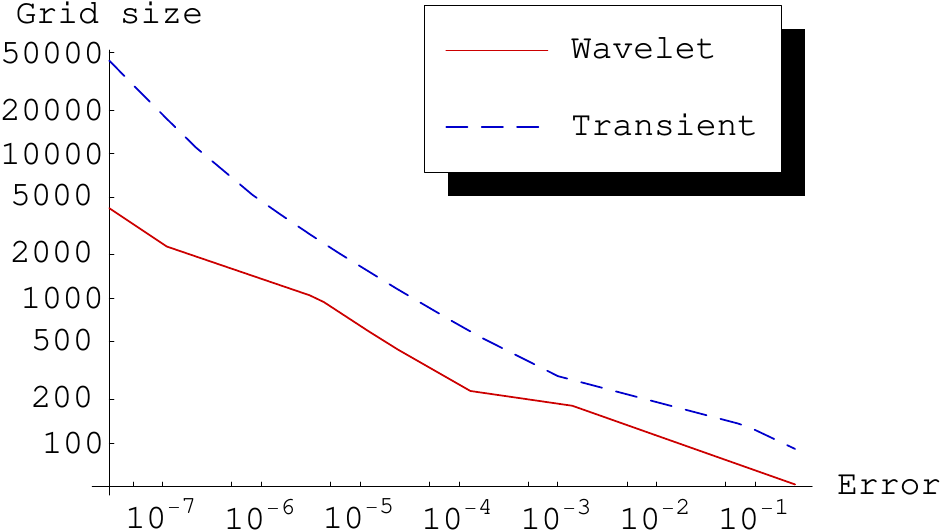}
\end{tabular}
\caption{Simulation results for the Schmitt trigger.
Computation time versus error (left), and
grid size versus error (right) for transient analysis and adaptive
wavelet analysis.
\label{bittner:schmitt_res}}
\end{figure}

{\bf Inverter chain.}
A further test circuit was an inverter chain consisting of 9
inverters. Therefore, the output signal represents the 9-times inverted
digital input signal. However, we can observe a delay and a
modification in the transition between high and low signal due to
intrinsic properties of used technology. Similar to the hysteresis effect 
in the previous problem, the
output depends strongly on the earlier behaviour of the input
signal, which again requires the use of the interval splitting
wavelet method to obtain the correct results.

\begin{figure}[htbp]
\centering
\begin{tabular}{c}
\includegraphics[width=0.9\textwidth]{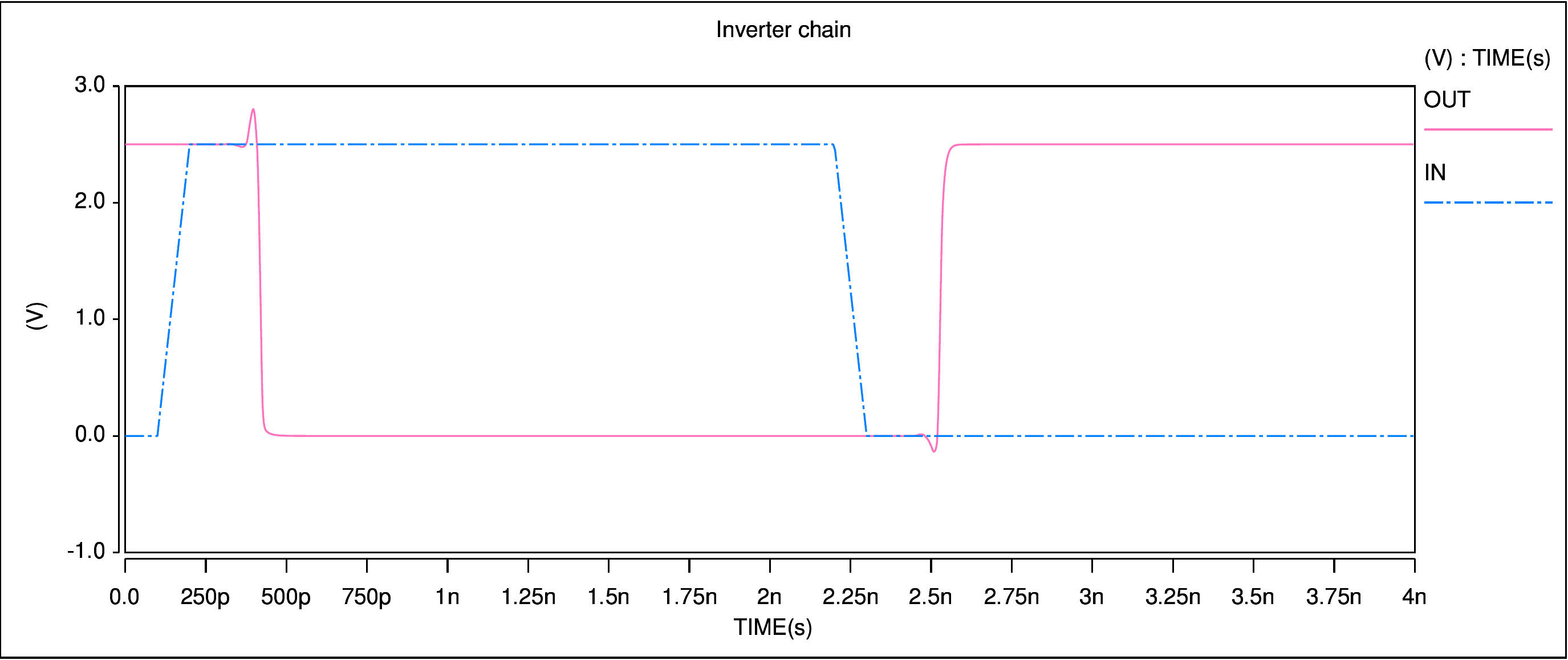}
\end{tabular}
\caption{Input and output signal of the inverter chain.
\label{bittner:inverter}}
\end{figure}

\begin{figure}[htbp]
\centering
\begin{tabular}{cc}
\includegraphics[width=0.48\textwidth]{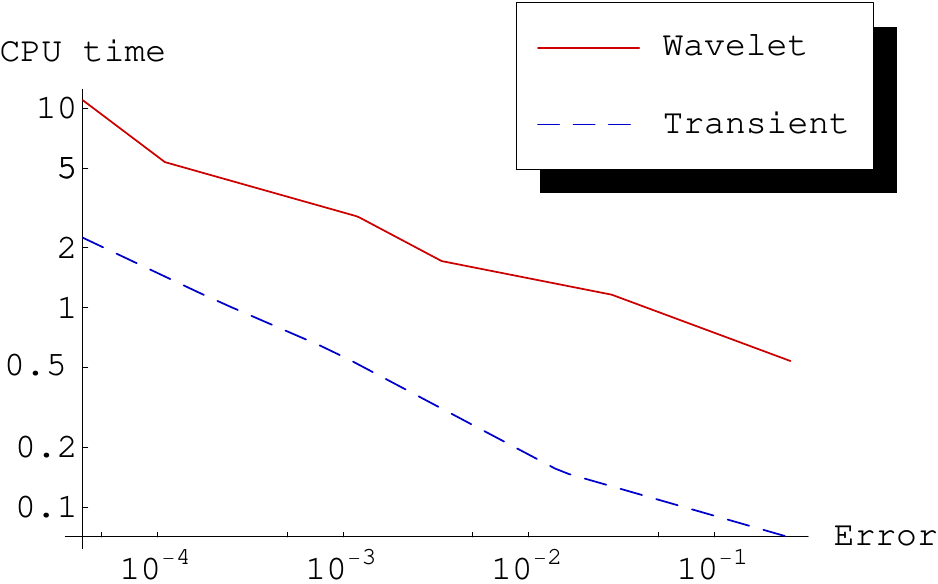}
 &\includegraphics[width=0.48\textwidth]{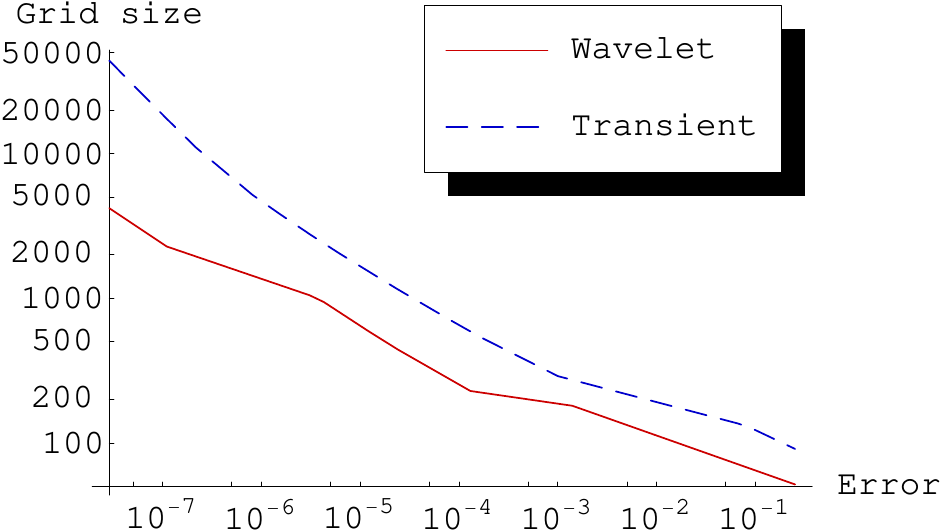}
\end{tabular}
\caption{Simulation results for the inverter chain.
Computation time versus error (left), and
grid size versus error (right) for transient analysis and adaptive
wavelet analysis.
\label{bittner:inv_res}}
\end{figure}

In both examples, the interval splitting wavelet method could
produce the correct results. i.e., we achieve robustness for the
wavelet-based approach. The performance is comparable to transient
analysis, although the current implementation is not faster than the
reference method. However, we see a big potential for the
improvement of the implemented method.

\section{Conclusion}

The results of the simulations indicate that the wavelet-based
method is able to fulfill all accuracy requirements and may achieve
the performance of the standard transient analysis. Since the
relatively new wavelet approach has a large potential for
optimization, we are optimistic that wavelet analysis will be a
valuable tool for circuit simulation in the future. Therefore our
activities on optimization and further development of the
wavelet-based algorithm are continuing.

\begin{acknowledgement}
  This work has been supported within the EU Seventh Research Framework
  Project (FP7) \mbox{ICESTARS} with with project number FP7/2008/ICT/214911.
\end{acknowledgement}

%
%




\end{document}